\documentclass[leqno, myheadings, twoside]{amsart}
\usepackage{cases}
\usepackage{amsmath, amsthm, amssymb, amscd, amsxtra,graphicx}
\usepackage{latexsym, amsfonts}
\usepackage{url}
\usepackage{texdraw}
\usepackage{epsfig}
\def\R{\mathbb{R}}

\def\H{\mathbb{H}}

\def\1-1{(1.1)}

\makeatletter \@addtoreset{equation}{section}

\makeatletter \renewcommand{\@biblabel}[1]{#1.}

\theoremstyle{remark}

\begin{document}
\title [Sharp higher-order Sobolev
inequalities in the hyperbolic space ${\H}^n$] {Sharp higher-order
Sobolev inequalities in the hyperbolic space ${\H}^n$}
\author{Genqian Liu}


\maketitle Department of Mathematics, Beijing Institute of
Technology,
 Beijing 100081, People's Republic of China.
 \ \
E-mail address:  liugqz@bit.edu.cn

\vskip 0.46 true cm

\vskip 15 true cm

\begin{abstract} \
 In this paper, we obtain the sharp $k$-th order Sobolev inequalities in the
  hyperbolic space ${\H}^n$ for all $k=1,2,3,\cdots$.
 This gives an answer to an open question raised by Aubin
 in [Aubin, Princeton University Press, Princeton (1982), pp.$\,$176-177] for $W^{k,2}({\H}^n)$
 with $k>1$. In addition, we prove that the associated Sobolev constants are optimal.
\end{abstract}

\vskip 1.3 true cm

\section{ Introduction}

\vskip 0.43 true cm

 In the Euclidean space ${\R}^n$, Aubin \cite{A3} and Talenti \cite{T}
  proved that for any $u\in C_0^\infty ({\R}^n)$,
 \begin{equation}\label{o-1}  \bigg(\int_{{\R}^n} |u(x)|^{2n/(n-2)} dx \bigg)^{(n-2)/n} \le
 \Lambda_1 \int_{{\R}^n} |\nabla u|^2 dx, \end{equation}
where the best Sobolev constant is \begin{eqnarray}\label{1-2}
\Lambda_1=\frac{4\omega_n^{-2/n}}{n(n-2) },\end{eqnarray}  and
$\omega_n$ is the surface area of the sphere  ${\mathbb{S}}^n=
\{x\in {\mathbb{R}}^{n+1}\big|
|x|=\sqrt{x_1^2+\cdots+x_{n+1}^2}=1\}$.
 The equality sign holds in (1.1) if and only if $u$ has the form:
 \begin{eqnarray*}   u(x)= \bigg[
\frac{2\epsilon}{\epsilon^2 +|x-x_0|^2}\bigg]^{\frac{n}{2}-1}, \;\;
\quad   x\in {\Bbb R}^n, \end{eqnarray*}
where $\epsilon >0$ and $x_0$ is any fixed point in ${\Bbb R}^n$.

In order to solve the famous Yamabe problem,
 Aubin established (see \cite{A1} and \cite{A2}) the sharp Sobolev inequality on the
sphere ${\mathbb{S}}^n$,  written as \begin{eqnarray}\label
{1-5}\quad \;\quad \; \bigg(\int_{{\mathbb{S}}^n}
|u(y)|^{2n/(n-2)} dy \bigg)^{(n-2)/n} \le
 \Lambda_1 \int_{{\mathbb{S}}^n} |\nabla u(y)|^2 dy + \omega_n^{-2/n}
 \int_{{\mathbb{S}}^n} |u(y)|^2dy \end{eqnarray}
 for any $u\in W^{1,2}({\mathbb{S}}^n)$.
 The equality sign holds if and only if $u$ has the form:
 \begin{eqnarray*} u(y)= \bigg[  \frac{1-\epsilon^2}{2\epsilon} \left(\frac{1+\epsilon^2}{1-\epsilon^2} -\cos r\right)
\bigg]^{1-\frac{n}{2}}, \quad \; 0\le r < \pi, \end{eqnarray*}
where $r$ is the distance from $y$ to $y_0$ on ${\mathbb{S}}^n$,
and $y_0$ is any fixed point on ${\mathbb{S}}^n$.

 In the hyperbolic space $({\H}^n, h)$, one has the following
Sobolev inequality (see \cite{H1} and \cite{H2}): \ \   for any $u\in
C_0^\infty ({\H}^n)$,
   \begin{eqnarray} \quad \;\, \bigg(\int_{{\H}^n} |u|^{2n/(n-2)} dV_h \bigg)^{(n-2)/n} \le
 \Lambda_1 \int_{{\H}^n} |\nabla u |^2 dV_h - \omega_n^{-2/n}
 \int_{{\H}^n} |u|^2dV_h. \end{eqnarray}

\vskip 0.25 true cm

In [5, p.$\,$176-177], T. Aubin raised the following open question:
 for $k>1$, can one establish the sharp Sobolev inequalities that
  are similar to (1.1), (1.3) and (1.4) for
  $W^{k,2}({\R}^n)$, $W^{k, 2}({\mathbb{S}}^n)$ and
$W^{k,2}({\H}^n)$, respectively?

\vskip 0.2 true cm

 Cotsiolis and Tavoularis \cite{CT} obtained
 the sharp higher-order Sobolev inequalities in ${\R}^n$:
  for $n>2k$ and any $u\in C_0^\infty({\R}^n)$,
     \begin{eqnarray} \label{1.5} \left(\int_{{\R}^n} |u(x)|^{\frac{2n}{n-2k}}
dx\right)^{\frac{n-2k}{n}} \le  \Lambda_k \int_{{\R}^n}
|\triangle^{k/2} u|^2 dx,\end{eqnarray} where
$$  |\triangle^{k/2} u|^2
 =\left\{
\begin{array}{ll}
 |\triangle^{k/2}u|^2   &\quad \;\; \mbox{if $
 k$ is even},\\
  |\nabla(\triangle^{(k-1)/2}u)|^2    &\quad \;\; \mbox{if $k$
 is odd}, \end{array}
\right. $$  and the best Sobolev constant is
\begin{eqnarray}\Lambda_k= \frac{2^{2k} \,\omega_n^{-(2k)/n}}{
  n\big[n-2k\big]\big[n^2 -\big(2(k-1)\big)^2\big]
  \big[n^2 -\big(2(k-2)\big)^2\big]
\cdots \big[n^2 - 2^2\big]}. \end{eqnarray} The equality sign holds
in (\ref{1.5}) if and only if
 \begin{eqnarray*} \label {1.7} u(x)=
 \bigg[
\frac{2\epsilon}{\epsilon^2 +|x-x_0|^2}\bigg]^{\frac{n}{2} -k}, \;\;
\quad   x\in {\Bbb R}^n. \end{eqnarray*}

 In \cite{B}, Beckner established the sharp $k$-th order Sobolev
 inequality on the sphere $({\mathbb{S}}^n, g)$:
  for $n>2k$  and any $u\in W^{k,2}({\mathbb{S}}^n)$,
\begin{eqnarray}\label{0.0}\left(\int_{{\mathbb{S}}^n} |u|^{\frac{2n}{n-2k}} dV_g
\right)^{\frac{n-2k}{n}} \le
   \Lambda_k \int_{{\mathbb{S}}^n} (Q_{k} u) u\, dV_g, \end{eqnarray}
where $Q_k$ is the  GJMS operator on ${\Bbb S}^n$:
$$Q_{k}=\frac{\Gamma\big(M+\frac{1}{2}+k\big)}{\Gamma
\big(M+\frac{1}{2}-k\big)},\quad \, \, M=
\sqrt{\triangle_{{\mathbb{S}}^n}+ \big(\frac{n-1}{2}\big)^2}.$$ The
equality sign holds in (\ref{0.0}) if and only if
\begin{eqnarray*} u(y)= \bigg[  \frac{1-\epsilon^2}{2\epsilon}
\left(\frac{1+\epsilon^2}{1-\epsilon^2} -\cos r\right)
\bigg]^{k-\frac{n}{2}}, \quad \; 0\le r < \pi, \end{eqnarray*}
  In 1995, Branson \cite{Br} gave the second proof for the sharp
  inequality (\ref{0.0}) by applying a more general Lie-theoretic point of view.
  Other proofs can be found in
  \cite{H} and \cite{CT}.

Another interesting problem is to discuss the  best Sobolev
constants.  Hebey (see \cite{H1}, \cite{H2}), Djadli, Hebey and
Ledoux \cite{DHL} have given best Sobolev constants for the
 first-order and second-order Sobolev inequalities on Riemannian
manifolds.

\vskip 0.2 true cm

 In this paper, for any positive integer $k$,
 we obtain the following sharp $k$-th order Sobolev inequality in the hyperbolic
space ${\H}^n$ of constant sectional curvature $-1$:

\vskip 0.30 true cm

 \noindent {\bf Theorem 1.1.}  {\it Let $({\H}^n,
h)$ be the hyperbolic $n$-space, $n>2k$, and let $q=(2n)/(n-2k)$.
Then, for any $u\in C_0^\infty({\H}^n)$,
\begin{eqnarray}\label {1.9} \left(\int_{{\H}^n} |u|^q dV_h \right)^{2/q} \le
   \Lambda_k \int_{{\H}^n} (P_k u) u\, dV_h, \end{eqnarray}
 where $P_k$ is a $2k$-th order  operator on ${\H}^n$ (see
Section 2) given by $P_k =P_1 (P_1+2) \cdots (P_1+k(k-1))$ with
$P_1=\triangle_h -\frac{n(n-2)}{4}$,
 $\, \triangle_h =-\frac{1}{\sqrt{|h|}}\sum_{i,j=1}^n
   \frac{\partial}{\partial x_i} \left( \sqrt{|h|}\,h^{ij}
   \frac{\partial}{\partial x_j}\right)$,
   and $\Lambda_k$ is the best $k$-th order Sobolev constant in ${\R}^n$.
 Moreover, for any $\epsilon>0$, if
\begin{eqnarray} \label{1---10} \psi_{k,\epsilon} (r)=\bigg[\frac{1+\epsilon^2}{2\epsilon}
\left(\cosh r- \frac{ 1-\epsilon^2}{1+\epsilon^2}
\right) \bigg]^{k-\frac{n}{2}},
 \, \quad 0\le r<+\infty,  \end{eqnarray}
 then
 \begin{eqnarray}\label{1.10}  \lim_{\epsilon\to 0^+}
 \frac{\int_{{\H}^n} (P_k \psi_{k,\epsilon}(r)) \psi_{k,\epsilon}(r)\, dV_h}{
\left(\int_{{\H}^n} |\psi_{k,\epsilon}(r)|^q dV_h \right)^{2/q}}= \inf_{u\in
C_0^\infty({\H}^n)\setminus \{0\}} \,\frac{\int_{{\H}^n} (P_k u) u\,
dV_h}{ \left(\int_{{\H}^n} |u|^q dV_h \right)^{2/q}},
\end{eqnarray} and  \begin{eqnarray}\label{1.11} P_k
 \psi_{k,\epsilon} (r) = \frac{1}{\Lambda_k \omega_n^{2k/n}} \,
\big(\psi_{k,\epsilon}(r)\big)^{q-1}, \quad  0\le
r<+\infty,\end{eqnarray}  where $r$ is the distance from $y$ to $0$ on ${\H}^n$.}

\vskip 0.18 true cm

 This answers the Aubin question mentioned above.
 In addition, if  (\ref{1.9}) is re-written  as
\begin{eqnarray}   \left(\int_{{\H}^n} |u|^q dV_h \right)^{2/q} \le
   \Lambda_k \int_{{\H}^n} \bigg( |\triangle_h^{k/2} u|^2 +
 \sum_{m=0}^{k-1} a_{km} |\triangle_h^{m/2} u|^2\bigg) dV_h, \end{eqnarray}
 then we can further prove that the constants $\Lambda_k, \Lambda_k a_{k,k-1}, \cdots,
\Lambda_k a_{k0}$
  are optimal because they cannot be lowered
  (see Theorem 3.3),  where $a_{km}$ are the coefficients of $P_k$ (i.e.,
  $P_k =\Delta_h^k +\sum_{m=0}^{k-1} a_{km} \Delta_h^{m}$).

\vskip 0.3 true cm

   The main idea of this paper is as follows. By using the conformal
  map $\sigma: B_n\to {\Bbb H}^n$ (see (2.1) of Section 2), we first lift the
    extremal functions $G_{k,\epsilon} (x)= \left[
\frac{2\epsilon}{\epsilon^2 +|x|^2}\right]^{\frac{n}{2} -k}$
   of the sharp inequalities (1.5)
   to the functions $\psi_{k,\epsilon} (r)=
\left[\frac{1+\epsilon^2}{2\epsilon} \left(\cosh r- \frac{ 1-\epsilon^2}{1+\epsilon^2}
\right) \right]^{k-\frac{n}{2}}$ in ${\H}^n$ (see Section 2).
 This transform preserves the $L^{\frac{2n}{n-2k}}$ norm (i.e., $\int_{{\Bbb H}^n}
 |\psi_{k,\epsilon} (r)|^{\frac{2n}{n-2k}} dV_h= \int_{B_n} |G_{k,\epsilon} (x)|^{\frac{2n}{n-2k}}dx$).
  Next, we search for a $2k$-th order linear differential operator $P_k$
  such that $$\int_{{\Bbb H}^n} (P_k \psi_{k,\epsilon} (r))\psi_{k,\epsilon}(r)dV_h=
  \int_{B_n} (\Delta^k G_{k,\epsilon} (x))G_{k,\epsilon}(x)\,dx.$$  To this end, we seek out an operator $P_k$ satisfying a stronger requirement:
  $$P_k \psi_{k,\epsilon} (r) = \frac{1}{\Lambda_k \omega_n^{2k/n}} (\psi_{k,\epsilon} (r))^{\frac{n+2k}{n-2k}}
  \quad \; \mbox{for} \;\; 0\le r<+\infty.$$
    By some direct calculations, we obtain the explicit
    expression of $P_k$ (i.e., $P_k=P_1(P_1+2)\cdots (P_1+k(k-1))$ with $P_1=\Delta_h -\frac{n(n-2)}{4}$)
    and prove that $P_k$
 is a conformal covariant differential operator, i.e.,
 \begin{eqnarray*}   (P_k u)\circ \sigma
 =J_\sigma^{-\frac{n+2k}{2n}}\triangle^k \big[J_\sigma^{\frac{n-2k}{2n}} \big(u\circ
 \sigma\big)\big], \quad  \; \mbox{for all}\;\; u\in C_0^\infty ({\Bbb H}^n). \end{eqnarray*}
  Finally, we  shall prove the sharp higher-order Sobolev inequalities (\ref{1.9}) in ${\H}^n$
 and show that the constants $\Lambda_k, \Lambda_k a_{k,k-1}, \cdots, \Lambda_k a_{k0}$
  cannot be lowered (see Section 3).

 \vskip 1.3true cm

 \section{An new method getting higher-order GJMS operators in the hyperbolic space
${\H}^n$}

\vskip 0.48 true cm

The hyperbolic $n$-space ${\H}^n$ ($n\ge 2$) is a complete simple
connected Riemannian manifold having constant sectional curvature
equal to $-1$, and  for a given dimensional number, any two such
spaces are isometric \cite{W}. There are several models for
${\H}^n$, the most important being the half-space model, the ball
model, and the hyperboloid or Lorentz model, with the ball model
being especially useful for questions involving rotational symmetry.
 We will only use the ball model in this paper.

Let $B_n=\{x=(x_1, \cdots,x_n)\in {\R}^n\big|
(x_1^2+\cdots+x_n^2)^{1/2}<1\}$ be the unit ball in the Euclidean
space ${\R}^n$. For $B_n$, if we endow with the Riemannian metric
$$ ds^2 :=\frac{4|dx|^2}{(1-|x|^2)^2},$$ then the sectional
curvature becomes the constant $-1$. Furthermore, if we now define
spherical coordinates about $x=0$ by
\begin{eqnarray}  x=t\zeta,
\quad \; t=\tanh (r/2),\end{eqnarray}
 where $t\in [0, 1)$,
 $r\in [0, \infty)$, $\zeta\in {\Bbb S}^{n-1}$, then we obtain
the metric
$$ ds^2= (dr)^2 +(\sinh^2 r) |d\zeta|^2.$$
Note that for each $x\in B_n$, $\, t$ and $r$ are the Euclidean and
the hyperbolic distances from $0$ to $x$, respectively.
 One easily sees that, in the ball model,  the geodesics emanating from the origin are
 given by straight lines emanating from the origin,
and their length to the boundary $\partial B_n$ is infinite.

Let $\triangle_h$ be the Laplacian on ${\H}^n$ with the metric tensor
$h$, and let $F\in C^2({\H}^n; {\Bbb R}^1)$ with
 $$F(y)= f(r,\zeta),$$
 where $y=\mbox{Exp}\, (r\zeta)$ is the exponent map.
 We have by direct calculation (see p.$\,$40 of \cite{Ch1}) that
$$(\triangle_h F) (y(r,\zeta)) = -(\sinh r)^{1-n} \frac{\partial}{\partial r} \bigg(
 (\sinh r)^{n-1} \frac{\partial f}{\partial r}\bigg)
   - (\sinh r)^{-2} {\mathcal{L}}_\zeta f,$$
 where, when writing ${\mathcal{L}}_\zeta f$, we mean that $f\big|{\text{S}(r)}$
is to be considered as a function on ${\mathbb{S}}^{n-1}$ with
associated Laplacian ${\mathcal{L}}$. If $f$ is a radial function on
$({\H}^n,h)$ (i.e., function that depends only on distance from $0$
on ${\H}^n$),  then the corresponding Laplacian takes the following
simple form (see also p.$\;$180-181] of \cite{Ch2}): \begin{eqnarray}
\triangle_h f(r,\zeta) =-(\sinh r)^{1-n} \frac{\partial}{\partial r}
 \left[(\sinh  r)^{n-1} \frac{\partial f}{\partial r}\right]. \end{eqnarray}
 Similarly, for any positive integer $m$ we can define the
$m^{\text{th}}$-iterated operator $\triangle_h^m$ on
 the set of radial functions as the following:
 for any radial function $f\in C^{2m}({\H}^n; {\Bbb R}^1)$,
 $$   \triangle_h^m f(r) =-(\sinh  r)^{1-n} \frac{\partial}{\partial r}
 \left[(\sin  r)^{n-1} \frac{\partial (\triangle_h^{m-1} f)}{\partial r}\right],
 \,\quad \; m=1,2,\cdots.$$

We know (see Section 1) that for any positive integer $k$,
\begin{eqnarray}  G_{k,\epsilon} (x)=
 \bigg[
\frac{2\epsilon}{\epsilon^2 +|x|^2}\bigg]^{\frac{n}{2} -k}, \;\;
\quad   x\in {\Bbb R}^n \end{eqnarray}
 are the extremal functions
for the sharp $k$-th order Sobolev inequality (1.5) in the Euclidean
space ${\R}^n$,
 and
\begin{eqnarray} \label{2??4}  \triangle^k \left(G_{k,\epsilon} (x)\right)
=\frac{1}{\Lambda_k\omega_n^{2k/n}} \big(G_{k,\epsilon} (x)\big)^{q-1} \quad \;
\text{in} \;\; {\R}^n, \end{eqnarray} where $q=\frac{2n}{n-2k}$.
However, there is not extremal function for the sharp $k$-th order
Sobolev inequality in the unit ball $B_n\subset {\R}^n$. More
  generally, we have the following

\vskip 0.35 true cm

\noindent {\bf Lemma 2.1.}  {\it  Let $\Omega\subset {\R}^n$ be a
bounded domain with smooth boundary ($n>2k$), and let
$q=\frac{2n}{n-2k}$. Assume that $\,\Xi_k(\Omega)$ is a constant defined by
\begin{eqnarray}   \label{2012} \frac{1}{\Xi_k(\Omega)}=\inf_{u\in W^{k,2}_0
(\Omega)\,\setminus \{0\}}\;\; \frac{\int_\Omega |\triangle^{k/2}
u|^2 dx}{(\int_\Omega |u|^q dx)^{2/q}}.\end{eqnarray} Then
$\Xi_k(\Omega)=\Lambda_k$. Moreover, there is not extremal
function in $W_0^{k,2}(\Omega)$} for (2.5).

\vskip 0.29 true cm

\noindent {\it Proof.} \ \  The  Sobolev imbedding theorem implies that
 the right-hand side of (\ref{2012}) is finite.
We shall show that this value is nonzero.
 If it is not this case, then there exists a sequence of functions $u_l\in W^{k,2}_0(\Omega)$
such that $\int_{\Omega} |u_l(x)|^q dx=1$ and $\int_\Omega |\Delta^{k/2} u_l|^2 dx
\to 0$ as $l\to +\infty$. It follows from the well known fact
(see \cite{Lio} or p.$\;$229 of \cite{Gi}) that
$\|\Delta^{k/2} u_l\|_{L^2 (\Omega)}^2 =\|u_l\|^2_{W^{k,2}_0(\Omega)}$, where
$\|u_l\|^2_{W^{k,2}_0 (\Omega)} =\sum_{|\alpha|=k}\|D^\alpha u_l\|_{L^2(\Omega)}^2$. By
 applying the Sobolev imbedding theorem again, we get  $\|u_l\|_{L^q(\Omega)}
\le C \|u_l\|_{W^{k,2}_0(\Omega)}$, where $C$ dependents only on
$\Omega$ and $k$. Therefore, $\|u_l\|_{L^q(\Omega)}\to 0$ as $l\to \infty$. This is a contradiction,
  and the claim is proved.

 Next, it is not difficult to see that if
$\Omega'$ is a translation of $\Omega$ in ${\R}^n$, then
$\Xi_k(\Omega)=\Xi_k(\Omega')$; if $\Omega'\subset \Omega$, then
$\frac{1}{\Xi_k(\Omega)}\le \frac{1}{\Xi_k(\Omega')}$ because of
$W^{k,2}_0(\Omega')\subset W^{k,2}_0(\Omega)$.
   By translation in ${\R}^n$,  we may assume that $B (0;s_1)
   \subset\Omega \subset B(0;s_2)$, where $B(0;s) =\{x\in {\R}^n,
   |x|<s\},\;\, 0<s_1<s_2$. Thus  $\frac{1}{\Xi_k(B(0;s_1))}\ge
   \frac{1}{\Xi_k(\Omega)}\ge \frac{1}{\Xi_k(B(0;s_2))}$. On the other
   hand, $\Xi_k(B(0;s))$ is independent of $s$. Indeed, let $u\in
  W^{k,2}_0 (B(0; s))$. By setting $\tilde{u}(x)= u(sx)$, we have
  $\tilde{u}\in W^{k,2}_0 (B_1)$ and
  $\frac{\int_{B(0;1)} |\triangle^{k/2} {\tilde{u}}|^2 dx} {\big(\int_{B(0;1)}
|{\tilde{u}}|^{q}dx\big)^{2/q}}= \frac{\int_{B(0;s)}
|\triangle^{k/2} u|^2 dx} {\big(\int_{B(0;s)}
|u|^{q}dx\big)^{2/q}}$. This implies $\frac{1}{\Xi_k(B(0;s))}=\frac{1}{\Xi_k(\Omega)}$
for all $s>0$.

  Now, we shall prove that $\Lambda_k= \Xi_k(\Omega)$. Since
  $\bigcup_{s>0} W^{k,2}_0(B(0;s))$ is dense in $W^{k,2}({\R}^n)$, there exists a
  sequence $\{u_j\}$ in $W^{k,2}_0(B(0; s_j))$ such that
$$\frac{\int_{B(0;s_j)} |\triangle^{k/2} u_j|^2 dx}
{\big(\int_{B(0;s_j)} |u_j|^{q}dx\big)^{2/q}}\to \frac{1}{\Lambda_k}
 \;\; \mbox{as} \,\, j\to \infty.$$
 From
 $\frac{\int_{B(0;s_j)} |\triangle^{k/2} u_j|^2 dx}
{\big(\int_{B(0;s_j)} |u_j|^{q}dx\big)^{2/q}}\ge
\frac{1}{\Xi_k(B(0;s_j))}=\frac{1}{\Xi_k(\Omega)}$, we get $\frac{1}{\Lambda_k} \ge
\frac{1}{\Xi_k(\Omega)}$. It is obvious that $\frac{1}{\Xi_k(\Omega)}\ge \frac{1}{\Lambda_k}$.
Hence $\Xi_k(\Omega)= \Lambda_k$.

 Finally, we show that there is not extremal function  in $W_0^{k,2}(\Omega)$ for (\ref{2012}).
 Suppose that there exists $u\in W^{k,2}_0
(\Omega)\setminus \{0\}$ such that  $\frac{\int_\Omega
|\triangle^{k/2} u|^2 dx} {\big(\int_\Omega |u|^{q}dx\big)^{2/q}}
=\frac{1}{\Xi_k}$. The extension of $u$
 by zero outside $\Omega$ attains a minimum of
 $\frac{\int_{{\R}^n} |\triangle^{k/2} u|^2 dx}{\big(\int_{{\R}^n}
|u|^{q}dx\big)^{2/q}}$. This is a contradiction since (2.3) is the
unique form of the extremal functions for the sharp Sobolev
inequalities (\ref{1.5}). \ \ \  $\qquad \qquad \quad \; \square$

\vskip 0.31 true cm

It is easy to verify that
\begin{eqnarray} \label{2.1} \quad\quad \;
\lim_{\epsilon\to 0^+} \frac{\int_{B_n} |\triangle^{k/2} G_{k,\epsilon}|^2
dx}{\big(\int_{B_n}
 |G_{k,\epsilon}|^{\frac{2n}{n-2k}}dx\big)^{\frac{n-2k}{n}}}=
 \inf_{u\in W_0^{k,2}(B_n)\setminus \{0\}}\;\; \frac{\int_{B_n}
 |\triangle^{k/2} u|^2 dx}{\big(\int_{B_n}
 |u|^{\frac{2n}{n-2k}}dx\big)^{\frac{n-2k}{n}}}=\frac{1}{\Lambda_k}.\end{eqnarray}

\vskip 0.15 true cm

Let $\sigma$ be the conformal map from the Euclidean ball $B_n$ to
the hyperbolic space ${\H}^n$ defined by (2.1). Then the Jacobian of
$\sigma$ is $J_\sigma(x)=\big(\frac{2}{1-|x|^2}\big)^n$. Suppose $f$ is
a smooth function defined in the Euclidean ball $B_n$. Lift it to
the hyperbolic space ${\H}^n$ by formula:
\begin{eqnarray}  (J_\sigma (x))^{\frac{n-2k}{2n}}    F(y) =  f(x), \quad x\in B_n, \;\,
y=\sigma(x)\in {\H}^n, \end{eqnarray} i.e.,
 \begin{eqnarray*} F(y)= \left(\frac{2}{1-|x|^2}\right)^{k-\frac{n}{2}}
 f(x), \quad x\in {B}_n, \;
\, y=\sigma(x)\in {\H}^n. \end{eqnarray*}
By this formula, we can lift every function from the Euclidean ball ${B}_n$ to ${\H}^n$.
 The main purpose of lifting a function $f$ (defined in $B_n$) to $F$ (defined in ${\H}^n$) is
   that this transform preserves $L^q$-norm, i.e.,
   \begin{eqnarray} \label{2?1} \int_{B_n} |f(x)|^q dx = \int_{{\H}^n} |F(y)|^q dV_h, \end{eqnarray}
 where $dV_h=J_\sigma(x) dx$ and $q=\frac{2n}{n-2k}$.

It is not difficult to check that the previous function $G_{k,\epsilon} (x)$
(defined in the Euclidean ball ${B}_n$)
 is lifted to the function  $\big[\frac{1+\epsilon^2}{2\epsilon} \big(\cosh r- \frac{ 1-\epsilon^2}{1+\epsilon^2}
\big) \big]^{k-\frac{n}{2}}$  (defined in ${\H}^n$). In fact, since
$|x|=\tanh \frac{r}{2}$, we have \begin{eqnarray*}
  \frac{1+|x|^2}{1-|x|^2}=
 \frac{1+\tanh^2 \frac{r}{2}}{1-\tanh^2 \frac{r}{2}}
 = \cosh r. \end{eqnarray*}
 It follows that for $\epsilon >0$,
\begin{eqnarray*} \quad\quad
\frac{1+\frac{1}{\epsilon^2}|x|^2}{1-|x|^2}=\frac{1+\epsilon^2}{2\epsilon^2}
\bigg(\frac{1+|x|^2}{1-|x|^2}-\frac{1-\epsilon^2}{1+\epsilon^2}  \bigg)= \frac{1+\epsilon^2 }{2\epsilon^2 }
 \left(\cosh r-\frac{1-\epsilon^2}{1+\epsilon^2}\right). \end{eqnarray*}
 Thus
 \begin{eqnarray} \label {2--9}  && \left(\frac{2}{1-|x|^2}\right)^{k-\frac{n}{2}}
G_{k,\epsilon} (x) =
 \bigg(\frac{2}{1-|x|^2}\bigg)^{k-\frac{n}{2}}
\bigg(\frac{2\epsilon}{\epsilon^2 + |x|^2} \bigg)^{\frac{n-2k}{2}} \\
 && \quad \quad  =  \bigg[ \frac{\epsilon \left(1+\frac{1}{\epsilon^2} |x|^2\right)}{1-|x|^2} \bigg]^{k-\frac{n}{2}}
=  \bigg[ \frac{1+\epsilon^2}{2\epsilon} \left( \cosh r - \frac{1-\epsilon^2}{1+\epsilon^2} \right) \bigg]^{k-\frac{n}{2}}
:= \psi_{k,\epsilon} (r),\nonumber \\
 && \quad \qquad \quad  \quad  \quad  x\in B_n,  \;\, 0\le r<+\infty.\nonumber  \end{eqnarray}

In what follows we shall look for a linear differential operator $P_k$ of order $2k$ defined on
 ${\H}^n$ such that \begin{eqnarray} \label{2?2}  \int_{B_n} G_{k,\epsilon}  (x) (\Delta^k G_{k,\epsilon}  (x)) \, dx =
 \int_{{\H}^n} \psi_{k,\epsilon} (r) (P_k \psi_{k,\epsilon}(r)) dV_h.\end{eqnarray}
 Actually, by virtue of (\ref{2??4}) and $\int_{B_n} |G_{k,\epsilon}(x)|^q dx=\int_{{\Bbb H}^n}|\psi_{k,\epsilon}(r)|^q dV_h$,
 it suffices for us to find  an operator $P_k$ such that
  \begin{eqnarray*}\label{2?3} P_k\psi_{k,\epsilon} (r)
=\frac{1}{\Lambda_k \omega_n^{2k/n}} \big(\psi_{k,\epsilon}(r)\big)^{q-1},
\quad \, 0\le r< +\infty. \end{eqnarray*}
 The following Lemma gives an explicit expression for the $P_k$:

\vskip 0.35 true cm

\noindent  {\bf Lemma 2.2.} \ \  {\it  Assume that $k$ is a positive
integer, and assume that $({\H}^n, h)$ is the hyperbolic space,
$\,n>2k$. Let
\begin{eqnarray*}\label{2-9}  \psi_{k,\epsilon}(r)=
\bigg[\frac{1+\epsilon^2}{2\epsilon} \left(\cosh r- \frac{ 1-\epsilon^2}{1+\epsilon^2}
\right) \bigg]^{k-\frac{n}{2}},
 \,  \quad 0\le  r < +\infty.\end{eqnarray*}
 If $\psi_{k,\epsilon} (r)$ satisfies the equation
 \begin{eqnarray}\label{2?4} P_k\psi_{k,\epsilon} (r)
=\frac{1}{\Lambda_k \omega_n^{2k/n}} \big(\psi_{k,\epsilon}(r)\big)^{\frac{n+2k}{n-2k}},
\quad \, 0\le r < +\infty, \end{eqnarray}
   then $P_k$ has the form:
 \begin{eqnarray} \label{2?5} P_{k}= P_1 \big[P_1
+ 2\big]\big[P_1+6\big] \cdots \big[P_1 + k(k-1)\big],\end{eqnarray}
 where $P_1=\triangle_h -\frac{n(n-2)}{4}$.}

\vskip 0.40 true cm

\noindent {\it Proof.} \ \   (i)  \ \   For $k=1$,  since $\psi_{1,\epsilon}(r) =
\left[\frac{1+\epsilon^2}{2\epsilon} \left(\cosh r- \frac{ 1-\epsilon^2}{1+\epsilon^2}
\right) \right]^{1-\frac{n}{2}}$ and
  $$\triangle_h \psi_{1,\epsilon} (r)  =-(\sinh r)^{1-n} \frac{\partial}{\partial r}
 \left[(\sinh  r)^{n-1} \frac{\partial \psi_{1,\epsilon} (r)}{\partial r}\right].$$
 It is easy to verify that
\begin{eqnarray*}   \left(\triangle_h  -\frac{n(n-2)}{4}\right)\psi_{1,\epsilon} (r)= \frac{n(n-2)}{4}
(\psi_{1,\epsilon} (r))^{\frac{n+2}{n-2}}.\end{eqnarray*}
 In fact, this result had been known in \cite{H1} and \cite{H2} (also see \cite{Cha}, \cite{LP}) since
 $\triangle_h  -\frac{n(n-2)}{4}$ is the Yamabe operator on $({\H}^n, h)$.

\vskip 0.36 true cm

 (ii) \ \  For $k=2$,  we can directly verify
    \begin{eqnarray*}  \left(\triangle_h -\frac{n(n-2)}{4}\right)\left(\triangle_h -\frac{n(n-2)}{4}+2\right)
    \psi_{2,\epsilon} (r)
     = \frac{n(n-4)(n^2-2^2)} {2^4}(\psi_{2,\epsilon}(r))^{\frac{n+4}{n-4}},\end{eqnarray*}
    where $\psi_{2,\epsilon}(r)=\big[\frac{1+\epsilon^2}{2\epsilon} \big(\cosh r- \frac{ 1-\epsilon^2}{1+\epsilon^2}
\big) \big]^{2-\frac{n}{2}}$.
    Actually, this had been verified in \cite{DHL}
since  $(\triangle_h -\frac{n(n-2)}{4})(\triangle_h -\frac{n(n-2)}{4}+2)$ is just
 the Paneitz-Branson operator on $({\H}^n, h)$ (The fourth order Paneitz-Branson operator $Z_g$
on a Riemannain manifold $(\mathcal{M}, g)$ is given by  (see \cite{DHL}, \cite{Pa} or \cite{Cha})
        \begin{eqnarray*}  Z_g u=\Delta_g^2 u -div_g\left(\frac{(n-2)^2+4}{2(n-1)(n-2)}
     R_g g -\frac{4}{n-2}\, Ric_g\right) du +
      \frac{n-4}{2}Q_g u,  \end{eqnarray*}
     where \begin{eqnarray*}  Q_g = \frac{1}{2(n-1)} \Delta_g R_g + \frac{n^3
     -4n^2 +16 n -16}{8(n-1)^2 (n-2)^2}
     R_g^2 -\frac{2}{(n-2)^2} |Ric_g|^2, \end{eqnarray*}
      and $R_g$ and $Ric_g$ denote respectively the scalar curvature and Ricci curvature of the metric $g$.)

 \vskip 0.33 true cm

 (iii) \ \  For $k=3$ and $\psi_{3,\epsilon}(r)=
  \left[\frac{1+\epsilon^2}{2\epsilon} \left(\cosh r- \frac{ 1-\epsilon^2}{1+\epsilon^2}
\right) \right]^{3-\frac{n}{2}}$, by a straightforward calculation we get
 \begin{eqnarray*} && \left(\triangle_h -\frac{n(n-2)}{4}\right)\left(\triangle_h -\frac{n(n-2)}{4}+2\right)
 \left(\triangle_h -\frac{n(n-2)}{4}+6\right) \psi_{3,\epsilon} (r) \\
    && \quad \quad \quad =  \frac{n(n-6)(n^2-4^2)(n^2-2^2)}
 {2^6}(\psi_{3,\epsilon}(r))^{\frac{n+6}{n-6}}.\end{eqnarray*}
 (Note that the sixth-order operator $P_3$ had also been obtained by Branson \cite{Br2} and W\"{u}nsch \cite{Wu}).

\vskip 0.3 true cm

 (iv)  \ \  We shall prove by induction that $P_k$ has the form  (\ref{2?5}),
 and the equation (\ref{2?4}) holds for all $k\ge 1$.
Indeed, suppose that
  for all $1\le j\le k$:
\begin{eqnarray} \label {2q13} P_1 \big[P_1 + 2\big]
\cdots \big[P_1 + j(j-1)\big]\psi_{j,\epsilon} (r) = \frac{1}{\Lambda_j \omega_n^{2j/n}}
 \big(\psi_{j,\epsilon} (r)\big)^{\frac{n+2j}{n-2j}}.
\end{eqnarray}
Then \begin{eqnarray*} \triangle_h \psi_{k+1,\epsilon} &=&
      -(\sinh \, r)^{1-n} \frac{\partial}{\partial r}
    \left\{ \big(\sinh\, r)^{n-1} \frac{\partial}{\partial
    r}\left[ \bigg(\frac{1+\epsilon^2}{2\epsilon} \big(\cosh r- \frac{ 1-\epsilon^2}{1+\epsilon^2}
\big) \bigg)^{k+1-\frac{n}{2}} \right]\right\}\\
   &=&   - \left[ -\frac{n(n-2)}{4}+k(k+1) \right]
\psi_{k+1,\epsilon} (r) - 2  \left(\frac{1-\epsilon^2}{2\epsilon}\right)
   k(k-\frac{n}{2}+1) \psi_{k,\epsilon}(r)
\\  & &+(k-\frac{n}{2}
   )(k-\frac{n}{2}+1)\psi_{k-1,\epsilon}(r).
\end{eqnarray*}
  Acting  on both sides of the  equation above by the $P_k$, we find  that
\begin{eqnarray*}  &&\left[\triangle_h - \frac{n(n-2)}{4}
+k(k+1)\right]P_k \psi_{k+1,\epsilon}(r)  =  - 2\left(\frac{1-\epsilon^2}{2\epsilon}\right)
  k(k-\frac{n}{2}+1)
  P_k \psi_{k,\epsilon} (r) \\
  &&  + (k-\frac{n}{2})(k-\frac{n}{2}+1)
     \left(\triangle_h -\frac{n(n-2)}{4}+k(k-1)\right) P_{k-1}\psi_{k-1,\epsilon}(r).
  \end{eqnarray*}
   From assumption (\ref{2q13}), we get that
\begin{eqnarray*}  && \left[\triangle_h - \frac{n(n-2)}{4}
+k(k+1)\right]P_k \psi_{k+1,\epsilon} =
   -  \frac{2}{\Lambda_k \omega_n^{2k/n}}\left(\frac{1-\epsilon^2}{2\epsilon}\right)
   k(k-\frac{n}{2}+1) \big(\psi_{k, \epsilon} (r)\big)^{\frac{n+2k}{n-2k}}\quad \\
  \\  &&  +  \frac{1}{\Lambda_{k-1}\omega_n^{2(k-1)/n}} (k-\frac{n}{2})(k-\frac{n}{2}+1)
   \big(-\frac{n(n-2)}{4}+k(k-1) + \triangle_h \big)
\big(\psi_{k-1,\epsilon}(r)\big)^{\frac{n+2(k-1)}{n-2(k-1)}} \quad \end{eqnarray*}
  \begin{eqnarray*} &=& - \, \frac{2}{\Lambda_{k}\omega_n^{2k/n}}\left(\frac{1-\epsilon^2}{2\epsilon}\right)
   k(k-\frac{n}{2}+1)
\left(\frac{1+\epsilon^2}{2\epsilon}\big(\cosh \,r-
\frac{ 1-\epsilon^2}{1+\epsilon^2}
 \big)\right)^{-k-\frac{n}{2}} + \quad \quad \quad\;
  \\  &&  + \frac{1}{\Lambda_{k-1}\omega_n^{2(k-1)/n}}
  (k-\frac{n}{2})(k-\frac{n}{2}+1)
   \big(-\frac{n(n-2)}{4}+k(k-1)\big)\big(\psi_{k-1,\epsilon}(r)\big)^{\frac{n+2(k-1)}{n-2(k-1)}}
    \quad \quad \quad \;\\
 && + \frac{1}{\Lambda_{k-1}\omega_n^{2(k-1)/n}}
  (k-\frac{n}{2})(k-\frac{n}{2}+1)
\left(\frac{1+\epsilon^2}{2\epsilon} \right)^{-k+1-\frac{n}{2}}\times \;\quad
\\ &&   \quad \; \times
  \left[ (\frac{n}{2}-k)(\frac{n}{2}+k-1)
 \left(\cosh \,r-
\frac{ 1-\epsilon^2}{1+\epsilon^2}
\right)^{-k+1-\frac{n}{2}}\right. \, \quad\; \quad \; \;\\
&& \quad \;\;\quad \left. -2\left(\frac{1-\epsilon^2}{1+\epsilon^2} \right)
k(\frac{n}{2}+k-1)\left(\cosh \,r-
\frac{ 1-\epsilon^2}{1+\epsilon^2}
\right)^{-k-\frac{n}{2}} \right. \quad \quad \;\;\\
 && \quad \quad \;\; \left. + \left(\frac{1+\epsilon^2}{2\epsilon}\right)^{-2}
 \big(\frac{n}{2}+k\big)\big(\frac{n}{2} +k-1\big) \left(\cosh \,r-
\frac{ 1-\epsilon^2}{1+\epsilon^2}
\right)^{-k-1-\frac{n}{2}}\right]\quad \quad \;\;
\\ &= &
  \frac{1}{\Lambda_{k-1}\omega_n^{2(k-1)/n}}\big(\frac{n}{2} -k+1)\big)
  \big(\frac{n}{2}+k\big)\big(\frac{n}{2}-k\big)\big(\frac{n}{2}+(k-1)\big)
  \big(\psi_{k+1,\epsilon}(r)\big)^{\frac{n+2(k+1)}{n-2(k+1)}} \quad \quad \;\;\\
 &= & \frac{1}{\Lambda_{k+1}\omega_n^{2(k+1)/n}}
 \big(\psi_{k+1,\epsilon}(r)\big)^{\frac{n+2(k+1)}{n-2(k+1)}}.\quad \quad \;\;
  \end{eqnarray*}
Thus, (\ref{2?4}) is true
 for the operator $P_k$ having the form (\ref{2?5}). $\qquad \quad \;\;\; \square$

\vskip 0.35 true cm

By applying  the identity (\ref{2??4}) and the lifting formula
$\psi_{k,\epsilon} (r)=\big(J_\sigma(x)\big)^{\frac{2k-n}{2n}} G_{k,\epsilon}  (x)$
 we immediately see that
 the equation (\ref{2?4}) is equivalent to the following equation:
 \begin{eqnarray} \label{2??6}  P_k\psi_{k,\epsilon} (r)= \big(J_\sigma (x)\big)^{-\frac{n+2k}{2n}}
 \Delta^k G_{k,\epsilon}  (x), \quad  \;  x\in B_n,\;\;  0\le r< +\infty. \end{eqnarray}

\vskip 0.15 true cm

 In general, we can  prove that (\ref{2??6}) also holds for all
$u\in W_0^{k,2}({\H}^n)$.

\vskip 0.30 true cm

 \noindent {\bf Theorem 2.3.}  \ \  {\it  Let $({\H}^n, h)$ be the hyperbolic
$n$-space with  $$ h :=\left(\frac{2}{1-|x|^2}\right)^2\delta,$$
where $\delta$ is the Euclidean metric of ${\R}^n$.
Let \begin{eqnarray}\label{2-12}
 P_{k}= P_1 \big[P_1 + 2\big]\big[P_1+6\big]
\cdots \big[P_1 + k(k-1)\big], \end{eqnarray}
where $P_1 =\triangle_h -\frac{n(n-2)}{4}$, and  $\triangle_h$ is the
Laplacian on  ${\H}^n$ which has the local representation:
 $$\triangle_h =-\frac{1}{\sqrt{h}} \sum_{i,j=1}^n \frac{\partial}{\partial x_i}
 \left(\sqrt{h}\, h^{ij} \frac{\partial }{\partial x_j}\right).$$
Suppose $\sigma: B_n\to {\H}^n$ is the conformal map
defined by (2.1). Then for any $u\in
C_0^\infty({\H}^n)$,
\begin{eqnarray} \label{2.15}  (P_k u)\circ \sigma
 =J_\sigma^{-\frac{n+2k}{2n}}\triangle^k \big[J_\sigma^{\frac{n-2k}{2n}} \big(u\circ
 \sigma\big)\big], \quad  \text{for} \;\;  x\in B_n, \end{eqnarray}
 where $J_\sigma(x)=\left(\frac{2}{1-|x|^2}\right)^n$ is the Jacobian of $\sigma$,
 and $\triangle^k$
is the standard $k^{\text{th}}$-iterated Laplacian in ${\R}^n$.}

\vskip 0.38 true cm

\noindent {\it Proof.}  \   We shall prove this theorem by induction.
 For the sake of convenience,
 we simply write $u\circ \sigma$ and $(P_k
 u)\circ \sigma$ as $u$ and $P_k u$ in $B_n$, respectively. Then, (\ref{2.15}) is re-expressed as
\begin{eqnarray} \label{2m1}   \xi_k^{\frac{n+2k}{n-2k}}  (P_k u)=
\triangle^k  \big(\xi_k u \big), \quad
\,\,\text{for}\;\; x\in B_n, \end{eqnarray}
 where $\xi_k (x)= \big(J_\sigma (x)\big)^{\frac{n-2k}{2n}}$.

For $k=1$, noting that $$\sqrt{h} =\bigg(\frac{2}{1-|x|^2}\bigg)^n,$$
 we have
  \begin{eqnarray*}  \triangle_h u  &= & -\sum_{i,j=1}^n \left[h^{ij}
  \frac{\partial^2 u} {\partial x_i \partial x_j}
+\frac{1}{\sqrt{h}}
 \frac{\partial u}{\partial x_j} \frac{\partial}{\partial x_i}
\left(\sqrt{h}\, {h}^{ij}\right)\right]\\
 &=& -\sum_{i,j=1}^n  \left[ \big(\frac{1-|x|^2}{2}\big)^2  \frac{\partial^2 u}{\partial x_i
\partial x_j} \delta_{ij} +
  \big(\frac{1-|x|^2}{2}\big)^n \frac{\partial u}{\partial x_j}\,
\frac{\partial}{\partial x_i}\bigg( \big(\frac{2}{1-|x|^2}\big)^{n-2} \delta_{ij}\bigg)\right]\\
 &=& \left( \frac{1-|x|^2}{2}\right)^2 \triangle u  - (n-2) \left(\frac{1-|x|^2}{2}\right) \sum_{i=1}^n
  x_i \frac{\partial u}{\partial x_i},\end{eqnarray*}
 so that
  \begin{eqnarray*} & &\xi_1^{\frac{n+2}{n-2}} (P_1 u)=
\left( \frac{2}{1-|x|^2}\right)^{\frac{n+2}{2}}
  \left(\triangle_{h}  u -\frac{n(n-2)}{4} u\right) \\
 &=&
\left(\frac{2}{1-|x|^2}\right)^{\frac{n+2}{2}}
\left[ \big(\frac{1-|x|^2}{2}\big)^2  \triangle u  - (n-2) \big(\frac{1-|x|^2}{2}\big)
\sum_{i=1}^n   x_i \frac{\partial u}{\partial x_i} -\frac{n(n-2)}{4}u\right]\\
   &=& \big(\frac{2}{1-|x|^2}\big)^{\frac{n-2}{2}} \triangle u
 - (n-2) \big( \frac{2}{1-|x|^2}\big)^{\frac{n}{2}} \sum_{i=1}^n  x_i \frac{\partial u}{\partial x_i}
- \frac{n(n-2)}{4} \big(\frac{2}{1-|x|^2}\big)^{\frac{n+2}{2}} u.
\end{eqnarray*}
  On the other hand, we get 
 $$\frac{\partial \xi_1}{\partial x_i} =   \left(\frac{n}{2}-1\right)
\bigg(\frac{2}{1-|x|^2}\bigg)^{\frac{n}{2}} x_i, \quad \; i=1, 2, \cdots, n$$
 and
 $$\frac{\partial^2 \xi_1}{\partial x_i^2} =
\bigg(\frac{n}{2}-1\bigg) \bigg(\frac{2}{1-|x|^2}\bigg)^{\frac{n}{2}}
 + \bigg(\frac{n}{2} -1\bigg)\bigg(\frac{n}{2}\bigg)
\bigg(\frac{2}{1-|x|^2}\bigg)^{\frac{n+2}{2}} x_i^2,$$
so that
$$\triangle \xi_1 = -\sum_{i=1}^n \frac{\partial^2 \xi_1}{\partial x_i^2}
 =-\frac{n(n-2)}{4} \big(\frac{2}{1-|x|^2}\big)^{\frac{n+2}{2}}.$$
  It follows that
\begin{eqnarray*}  \triangle (u \xi_1)&=& -\sum_{i=1}^n
\frac{\partial^2 (u \xi_1)}{\partial x_i^2}
= \xi_1\triangle u -2 \nabla u\cdot \nabla \xi_1 +u\triangle \xi_1\\
    &= & \big(\frac{2}{1-|x|^2}\big)^{\frac{n-2}{2}} \triangle u
 - (n-2) \big( \frac{2}{1-|x|^2}\big)^{\frac{n}{2}} \sum_{i=1}^n
 x_i \frac{\partial u}{\partial x_i}\\
   && - \frac{n(n-2)}{4} \big(\frac{2}{1-|x|^2}\big)^{\frac{n+2}{2}} u.\end{eqnarray*}
 Thus  \begin{eqnarray} \label{2a0} \bigg(\triangle_h -\frac{n(n-2)}{4}\bigg)u
 =\bigg(\frac{2}{1-|x|^2}\bigg)^{-\frac{n+2}{2}}
 \triangle\bigg(\big(\frac{2}{1-|x|^2}\big)^{\frac{n-2}{2}}u\bigg),\end{eqnarray}
   i.e., (\ref{2m1}) is true  for $k=1$.

\vskip 0.25 true cm

   We now suppose that (\ref{2m1}) is true for $k$, that is,
    \begin{eqnarray}\label{2-18} \left( \frac{2}{1-|x|^2}
 \right)^{\frac{n+2k}{2}} (P_k u)= \triangle^k \left(\bigg(
 \frac{2}{1-|x|^2}\bigg)^{\frac{n}{2}-k} u\right), \end{eqnarray}
  and attempt to deduce from this that (\ref{2m1}) is still true for $k+1$.
   By a simple calculation, \ we have  that
  \begin{eqnarray} \label{2a1}
  && \triangle \left[\left(\frac{1-|x|^2}{2}\right)^{k+1} \triangle^k
 v \right]= (k+1)(n+2k)\bigg(\frac{1-|x|^2}{2}\bigg)^k \triangle^k
 v \\ && \quad \;\; \, + \left(\frac{1-|x|^2}{2}\right)^{k+1} \triangle^{k+1}
 v
 + 2(k+1)\left(\frac{1-|x|^2}{2}\right)^{k} \sum_{i=1}^n x_i \frac{\partial (\triangle^k
 v)}{\partial x_i} \nonumber \\
 && \quad \;\; \, -k(k+1)\left(\frac{1-|x|^2}{2}\right)^{k-1} \triangle^k
 v. \nonumber \end{eqnarray}
 It is easy to check that for each integer $k\ge 0$,
\begin{eqnarray} \label{2a2}
  && \triangle^{k+1} \left(\frac{1-|x|^2}{2} v\right)= (k+1)(n+2k)
   \triangle^k  v  \\
 && \quad \quad \; + \left(\frac{1-|x|^2}{2}\right) \triangle^{k+1}
 v  +2(k+1) \sum_{i=1}^n x_i \frac{\partial (\triangle^k
 v)}{\partial x_i}. \nonumber \end{eqnarray}
Combining (\ref{2a1}) and (\ref{2a2}), we get
\begin{eqnarray} \label{2a3}
  && \triangle \left[ \left(\frac{1-|x|^2}{2}\right)^{k+1} \triangle^k
 v \right] +k(k+1)\left(\frac{1-|x|^2}{2}\right)^{k-1} \triangle^k
 v \\
  && \quad \;\;  =  \left(\frac{1-|x|^2}{2}\right)^{k} \triangle^{k+1}
 \left(\frac{1-|x|^2}{2}v\right). \nonumber \end{eqnarray}
    By (\ref{2-12}), (\ref{2-18}), (\ref{2a0})  and (\ref{2a3}), we have that
 \begin{eqnarray*}
P_{k+1}u &=& \big[P_1 +k(k+1)\big]\, P_k u\\
&=& \bigg (\triangle_h -\frac{n(n-2)}{4}\bigg) P_k u +k(k+1)P_k u\\
 &=&\bigg(\triangle_h -\frac{n(n-2)}{4}\bigg)
 \left[ \left(\frac{2}{1-|x|^2}\right)^{-\frac{n+2k}{2}}
\triangle^k
\left(\bigg(\frac{2}{1-|x|^2}\bigg)^{\frac{n-2k}{2}}u\right)\right]\\
&& + k(k+1) \left(\frac{2}{1-|x|^2}\right)^{-\frac{n+2k}{2}}
\triangle^k
\left(\bigg(\frac{2}{1-|x|^2}\bigg)^{\frac{n-2k}{2}}u\right)\\
  &=&
\bigg(\frac{2}{1-|x|^2}\bigg)^{-\frac{n+2}{2}} \triangle
\left[\bigg(\frac{2}{1-|x|^2}\bigg)^{-(k+1)} \triangle^k
\left(\bigg(\frac{2}{1-|x|^2}\bigg)^{\frac{n-2k}{2}}u\right)\right]\\
&&  + k(k+1) \left(\frac{2}{1-|x|^2}\right)^{-\frac{n+2k}{2}}
\triangle^k
\left(\bigg(\frac{2}{1-|x|^2}\bigg)^{\frac{n-2k}{2}}u\right)\\
  &=&
\bigg(\frac{2}{1-|x|^2}\bigg)^{-\frac{n+2}{2}}
\bigg(\frac{2}{1-|x|^2}\bigg)^{-k} \triangle^{k+1}
\left(\bigg(\frac{2}{1-|x|^2}\bigg)^{\frac{n-2(k+1)}{2}}
u\right)\\
&=& \bigg(\frac{2}{1-|x|^2}\bigg)^{-\frac{n+2(k+1)}{2}}
\triangle^{k+1}
\left(\bigg(\frac{2}{1-|x|^2}\bigg)^{\frac{n-2(k+1)}{2}} u\right).
\end{eqnarray*}
 Hence $$\xi_{k+1}^{\frac{n+2(k+1)}{n-2(k+1)}}(P_{k+1}u)
 =\triangle^{k+1} (\xi_{k+1}\,u),$$
 which  completes the proof. $\,\, \square$

\vskip 0.55 true cm

\noindent {\bf Remark 2.4}. \  \  (i) \  \ {\it  In the general setting of psuedo-Riemannian
manifolds, Graham,
  Jenne, Mason and Sparling solved a major existence problem in
  \cite{GJMS} where they used a formal geometric construction to show
  the existence of conformally covariant differential operators $P_k$ (to be
  referred to as the GJMS operators) with principal part $\triangle^k$.

(ii) \ \  In fact, similar explicit representation
 of the GJMS operator $P_k$ can be found for Einstein manifolds
 in \cite{FG}, \cite {G} or \cite
 {GS}  by using the GJMS operator construction and some
 special properties of Einstein metrics.
However our approach for obtaining $P_k$ in the hyperbolic space ${\H}^n$ is
 completely different. In this paper, we adopt an elementary
 calculation by a lifting function $\psi_{k,\epsilon}(r)$
  and by using the nonlinear elliptic equation (\ref{2?4})
which is equivalent to the equation (\ref{2??6}).
 By our new method, we can clearly see that the $2k$-th order GJMS operators
 $P_k$ on the hyperbolic space $({\H}^n, h)$ is essentially the power
 $\Delta^k$ of the positive ${\Bbb R}^n$ Laplacian lifted to the hyperbolic space via the
 conformal map $\sigma$.}

\vskip 0.31 true cm

 \vskip 1.49 true cm

\section{ Sharp $k$-th order Sobolev inequalities in the
hyperbolic space ${\H}^n$}

\vskip 0.48 true cm

 \noindent {\bf Lemma 3.1.} \ \ {\it  Let $({\H}^n, h)$ be the hyperbolic space,
$n>2k$.
  Suppose that there exist real constants
$b'_k$ and $\{a'_{km}\}_{0\le m\le k-1}$ such that for any $u\in
C_0^\infty ({\H}^n)$,
\begin{eqnarray}\label{3-1}  \quad\;\; b'_k
\left(\int_{{\H}^n} |u|^q dV_h\right)^{2/q} \le
\int_{{\H}^n}|\triangle_h^{k/2} u|^2dV_h  + \sum_{m=0}^{k-1}
      a'_{km} \int_{{\H}^n} |\triangle_h^{m/2} u|^2 dV_h, \end{eqnarray}
 where  $q=(2n)/(n-2k)$.
 Then $b'_k \le \frac{1}{\Lambda_k}$,
where $\Lambda_k$ is  the best $k$-th order Sobolev constant in
${\R}^n$.}

\vskip 0.28 true cm

\noindent {\it Proof.} \  This proof follows from the lines of
\cite{A1} (also see \cite{H1}).  Suppose by contradiction that there
exist
     $b'_k$ and  $\{a'_{km}\}$ satisfying $b'_k>\frac{1}{\Lambda_k}$, such that
 inequality (3.1) holds
 for any $u\in C_0^\infty({\H}^n)$.
 Let $y\in {\H}^n$. It is easy to see that for any $\epsilon>0$ there exists
a chart $(U, \phi)$ of ${\H}^n$ at $y$, and there exists $\tau>0$
such that $\phi (U)= B_\tau (0)$ (here $B_\tau (0)$ is the
Euclidean ball of center $0$ and radius $\tau$ in ${\R}^n$), and
such that the components $h_{ij}$ of $h$ in this chart satisfy
$$ (1-\epsilon)\delta_{ij} \le h_{ij} \le (1+\epsilon)\delta_{ij}$$
as bilinear forms.
   Choosing $\epsilon$ small enough we can get by (3.1)
 that there exist $\tau_0>0$, $\, b''_k$ and $\{a''_{km}\}$ satisfying $b''_{kk}
 > \frac{1}{\Lambda_k}$
   such that
for any $\tau\in (0,\tau_0)$ and any $u\in C_0^\infty(B_\tau (0))$,
 $$  b''_k \left(\int_{{\R}^n} |u|^q dx\right)^{2/q} \le
   \int_{{\R}^n} |\triangle_h^{k/2} u|^2 dV_h + \sum_{m=0}^{k-1}
    a''_{km}  \int_{{\R}^n}|\triangle^{m/2} u|^2 dx. $$
 Applying Nirenberg's lemma (see \cite{N} or [1, Lemma 14.1]), we find that
 there exists a constant $c$ depending only
 on $\tau$ such that for any $0< m \le k-1$,
$$ \int_{{\R}^n} |\triangle^{m/2} u|^2 \le
\tau \int_{{\R}^n} |\triangle^{k/2} u|^2 dx + c
 \int_{{\R}^n} |u|^2 dx.$$
It follows from H\"older's inequality that
$$\int_{B_\tau (0)} |u|^2 dx \le  |B_\tau(0)|^{2k/n} \left(\int_{B_\tau(0)}
   |u|^q\right)^{2/q},$$
where $|B_\tau(0)|$ denotes the volume of the ball $B_\tau(0)$ in
${\R}^n$. Therefore by choosing $\tau$ small enough, we obtain that
there exist $\tau>0$ and $b'''_k>\frac{1}{\Lambda_k}$ such that for
any $u\in C_0^\infty(B_\tau(0))$,
 $$ b'''_k \left(\int_{{\R}^n} |u|^q dx \right)^{2/q} \le
  \int_{{\R}^n} |\triangle^{k/2} u|^2dx.$$
For any $u\in C_0^\infty({\R}^n)$, let us set $u_\eta (x)=u(\eta x),
\; \eta >0$. Take $\eta$ large enough, $u_\eta\in C_0^\infty(B_\tau
(0))$.
 Thus,
$$  b'''_k \left(\int_{{\R}^n} |u_\eta|^q dx \right)^{2/q}\le
  \int_{{\R}^n} |\triangle^{k/2} u_\eta|^2dx. $$
  But we also have
$$\left (\int_{{\R}^n}|u_\eta|^q dx\right)^{2/q} =\eta^{-(2n)/q} \left(\int_{{\R}^n}
 |u|^q dx\right)^{2/q} $$
 and $$\int_{{\R}^n}|\triangle^{k/2} u_\eta|^2 dx =\eta^{2k-n} \int_{{\R}^n}
|\triangle^{k/2} u|^2 dx. $$
 In view of $1/q=1/2-k/n$, we obtain that for any $u\in C_0^\infty({\R}^n)$,
$$ b'''_k  \left(\int_{{\R}^n}|u|^qdx\right)^{1/q} \le  \int_{{\R}^n}
|\triangle^{k/2} u|^2 dx. $$ Since $b'''_k>\frac{1}{\Lambda_k}$,
such an inequality is contradiction with (1.5). This completes
proof of the lemma.$\;\;\square$

\vskip 0.43 true cm

\noindent  {\bf  Proof of theorem 1.1.} \ \
 We simply write  $u\circ \sigma$ and $P_k (u\circ \sigma)$
 as  $u$ and $P_k u$ in $B_n$, respectively, where $\sigma: B_n \to {\H}^n$ is the conformal map
  defined by (2.1).
   By Theorem 2.3 (or (2.17)) we have
\begin{eqnarray}\label{3-8} \xi_k^{\frac{n+2k}{n-2k}}(P_k u) =\triangle^{k} (\xi_k u)
\quad \, \text{for} \;\; x\in B_n, \end{eqnarray}
 where
$\xi_k=\left(\frac{2}{1-|x|^2}\right)^{\frac{n-2k}{2}}$, and
$\triangle^k$ is the $k^{\text{th}}$-iterated standard Laplacian in
${\R}^n$.
  Multiplying $\xi_k u$ to equation (\ref{3-8}) and then integrating the result in
 $B_n$, we get
 \begin{eqnarray}\label{3-9}  \int_{B_n} \left(\frac{2}{1-|x|^2}\right)^n (P_k u)u\, dx =
\int_{B_n} (\xi_k  u)\big(\triangle^k (\xi_k u)\big)
dx.\end{eqnarray}
By applying integration by parts to the right-hand side of (\ref{3-9}), we have
 \begin{eqnarray}\label{3-10}  \int_{B_n} (\xi_k  u)\big(\triangle^k ( \xi_k u)\big) dx=
 \int_{B_n} |\triangle^{k/2} (\xi_k u)|^2 dx.\end{eqnarray}
 Note that   $$  dV_h= \bigg(\frac{2}{1-|x|^2}\bigg)^n dx,$$
where $dx$ is the volume element of the Euclidean space ${\R}^n$.
This implies
 \begin{eqnarray}\label {3b11} \int_{B_n} \left(\frac{2}{1-|x|^2}\right)^n (P_k u)u\, dx
 =\int_{{\H}^n}  (P_k u)u\, dV_{h}. \end{eqnarray}
   It follows from (\ref{3-9})---(\ref{3b11}) that
\begin{eqnarray}\label{3-12}   \frac{\int_{{\H}^n} (P_k u)u\, dV_h}{\left(\int_{{\H}^n}
|u|^q dV_h\right)^{2/q}}=
  \frac{ \int_{B_n} |\triangle^{k/2} (\xi_ku)|^2 dx}
  {\left(\int_{B_n} |\xi_k \,u |^q dx\right)^{2/q}}.\end{eqnarray}
    It is clear that for every $u\in C_0^\infty({\H}^n)$,
\begin{eqnarray}\label{3-6} \frac{\int_{{\H}^n} (P_k u)u\, dV_h}{\left(\int_{{\H}^n}
|u|^q dV_h\right)^{2/q}}
  = \frac{\int_{{\H}^n} \left(|\triangle_{h}^{k/2} u|^2
+\sum_{m=0}^{k-1} a_{km}|\triangle_{h}^{m/2}
 u|^2\right)  dV_{h}}
{\left(\int_{{\H}^n} |u|^q dV_{h}\right)^{2/q}},\end{eqnarray}
 where $a_{km}$ are the coefficients of $P_k$.
Also, there exits a constant $C>0$ such that for any $u\in C_0^\infty({\Bbb H}^n)$
$$ \left(\int_{{\H}^n} |u|^q dV_h\right)^{2/q}\le C
\int_{{\H}^n} \left(|\triangle_h^{k/2} u|^2  +
\sum_{m=0}^{k-1} a_{km}  |\triangle_h^{m/2} u|^2\right) dV_h.$$
 From Lemma 3.1, we see that
 $$  \inf_{u\in C_0^\infty ({\H}^n)\setminus \{0\}}
 \frac{\int_{{\H}^n} \left(|\triangle_h^{k/2} u|^2  +
\sum_{m=0}^{k-1} a_{km}  |\triangle_h^{m/2} u|^2\right) dV_h }
{\left(\int_{{\H}^n} |u|^q dV_h\right)^{2/q}}\le
\frac{1}{\Lambda_k}.$$

Suppose by contradiction that
   $$  \inf_{ u\in C_0^\infty({\H}^n)\setminus \{0\}}
 \frac{\int_{{\H}^n} \left( |\triangle_h^{k/2} u|^2 +
\sum_{m=0}^{k-1} a_{km} |\triangle_h^{m/2} u|^2\right) dV_h }
{\left(\int_{{\H}^n} |u|^q dV_h\right)^{2/q}}<
\frac{1}{\Lambda_k},$$
 and let  $u_0\in  C_0^\infty({\H}^n)$, $u_0\not\equiv 0$ satisfy
 \begin{eqnarray}\label{3-13}   \frac{\int_{{\H}^n}  |\triangle_h^{k/2} u_0|^2 dV_h +
\sum_{m=0}^{k-1} a_{km} |\triangle_h^{m/2} u_0|^2 dV_h }
{\left(\int_{{\H}^n} |u_0|^q dV_h\right)^{2/q}}<
\frac{1}{\Lambda_k}.\end{eqnarray}
   By (\ref{3-12}), (\ref{3-6}) and (\ref{3-13}), we have
$$ \frac{\int_{B_n} |\triangle^{k/2} (\xi_k u_0)|^2  dx}
{ \left(\int_{B_n} |\xi_k u_0|^q
        dx\right)^{2/q}} <\frac{1}{\Lambda_k}.$$
 Clearly, $\xi_k u_0 \in C_0^\infty (B_n)$.
 But
$$  \frac{\int_{B_n} |\triangle^{k/2} (\xi_k u_0)|^2
dx} {\big(\int_{B_n}|\xi_k u_0|^q dx\big)^{2/q}} \ge \inf_{w\in
C_0^\infty(B_n)\setminus \{0\}} \;\frac{\int_{B_n} |\triangle^{k/2}
w|^2 dx} {\big(\int_{B_n}|w|^q dx\big)^{2/q}} \ge
\frac{1}{\Lambda_k}.$$
 This is a contradiction,
  which shows
  $$  \inf_{ u\in C_0^\infty({\H}^n)\setminus \{0\}}
 \frac{\int_{{\H}^n} \left( |\triangle_h^{k/2} u|^2 +
\sum_{m=0}^{k-1} a_{km} |\triangle_h^{m/2} u|^2\right) dV_h }
{\left(\int_{{\H}^n} |u|^q dV_h\right)^{2/q}}
  =\frac{1}{\Lambda_k}.$$
  Hence we get (\ref{1.9}).

Finally,  by the conformal map $\sigma$  from $(B_n, \delta)$ to
$({\H}^n, h)$, the function $G_{k,\epsilon}  (x)$ (see (2.3))
 defined in $B_n$ is lifted  to $\psi_{k,\epsilon} (r)$,
 and (\ref{1.10}) is true.
 (\ref{2?4}) of Lemma 2.2 also shows that equation (\ref{1.11}) holds.$\;\;\square$

\vskip 0.34 true cm

 \noindent {\bf Remark 3.2} \  \  (i) \  \  {\it Using a similar argument as in Section 2,
 we can conclude that there is not extremal
  function  on ${\H}^n$ for the sharp Sobolev inequality (1.8).

(ii)   \ \  If a Riemannian manifold $(\mathcal {M},g)$ is conformally flat (i.e., there exists a
 smooth function $w$ defined on $\mathcal{M}$ such that
 $(\mathcal{M}, e^{2w}g)$ is flat (i.e., the curvature of $e^{2w}g$ vanishes on $\mathcal{M}$)),
 then by a completely similar discussion as we have done in the hyperbolic space $({\H}^n, h)$,
  we can also obtain the corresponding sharp $k$-th order Sobolev inequalities on $(\mathcal{M}, g)$.
 In particular, we can obtain that if $(\mathcal{M},g)$ is a conformally flat
 Riemannian manifold and if $(\mathcal{M}, g)$ is also
 an Einstein manifold (i.e., $Ric_g =\beta g$  for some constant $\beta$),
 then the $2k$-th order GJMS operator $E_k$ has an explicit expression (also see \cite{G}):
\begin{eqnarray*} E_k= \prod_{l=1}^k (\Delta_g -c_l R_g), \end{eqnarray*}
where $c_l= \frac{(n+2l-2)(n-2l)}{4n(n-1)}$.
 The sharp $k$-th order Sobolev inequalities is
\begin{eqnarray*}
 \left(\int_{\mathcal{M}} |u|^q dV_g \right)^{2/q} \le
   \Lambda_k \int_{\mathcal{M}} (E_k u) u\, dV_g, \,\quad \forall \;
 u\in C_0^\infty(\mathcal{M}). \end{eqnarray*}
 More particularly, we can immediately obtain the GJMS operator $E_k$ on the unit sphere $({\Bbb S}^n,g)$
 (that is, $E_k =E_1(E_1-2)\cdots (E_1-k(k-1))$ with $E_1=\Delta_g +\frac{n(n-2)}{4}$, see
\cite{L} and \cite{Mo}).}

\vskip 0.39 true cm

  Recall that the inequality (\ref{1.9}) can be re-written as
 \begin{eqnarray*}   \|u\|_{L^q({\H}^n)}^2 \le \Lambda_k \left[\int_{{\Bbb H}^n} |\Delta^{k/2} u|^2 dV_h+
\sum_{m=0}^{k-1} a_{km} \int_{{\H}^n}
|\Delta_h^{m/2} u|^2 dV_h \right], \quad \,
  \forall \; u\in C_0^\infty({\H}^n), \end{eqnarray*}
where $a_{km}$  ($0\le m\le k-1$) are the coefficients of the operator $P_k$.
Theorem 1.1 only shows that $\Lambda_k$ is the best constant. The
following theorem says
  that the constants $\Lambda_k, \Lambda_k a_{k,k-1}, \cdots, \Lambda_k a_{k0}$
  are optimal since they cannot be lowered.

\vskip 0.28 true cm

\noindent  {\bf Theorem 3.3} \ \  {\it Let $({\H}^n, h)$ be the
hyperbolic space,
 $n>4k-2$,  and let $q=(2n)/(n-2k)$.
 Assume that $T_i$  ($i=0, 1, \cdots, k-1$) is the operator defined by
$$T_i =\sum_{m=i+1}^{k} a_{km} \triangle_h^m + \sum_{m=0}^{i} \tau_{km}\triangle_h^{m}, \quad \;\;
 i=0, 1,\cdots, k-1,$$
where $\{a_{km}\}$
 are the coefficients of the $2k$-th order GJMS operator $P_k$ with $a_{kk}=1$.
For each fixed $i\in \{0,1, \cdots, k-1\}$, there
exist real numbers $\tau_{k0}, \cdots, \tau_{ki}$ such that
   for all $u\in C_0^\infty({\H}^n)$,
 \begin{eqnarray} \label{3x1}  \|u\|_{L^q({\H}^n)}^2 &\le & \Lambda_k \int_{{\H}^n} (T_i  u)u\,
 dV_h \nonumber
\\  &= &\Lambda_k \left[\sum_{m=i+1}^k a_{km} \int_{{\H}^n}
|\triangle_h^{m/2} u|^2 dV_h +\sum_{m=0}^i \tau_{km} \int_{{\H}^n}
|\triangle_h^{m/2} u|^2dV_h \right] \end{eqnarray}
 if and only if $\tau_{ki}\ge a_{ki}$, where $\Lambda_k$ is the best $k$-th order Sobolev
constant in ${\R}^n$.}

\vskip 0.28 true cm

\noindent {\it Proof. }   For any fixed $i\in \{0, 1\cdots, k-1\}$,
if $\tau_{ki}\ge a_{ki}$, the result immediately follows from
  Theorem 1.1 because we may take $\tau_{km}=a_{km}$, $\, m=0,1,
\cdots, i-1$.

 Suppose, on the contrary, that there exist $\tau_{k0}, \cdots, \tau_{ki}$ satisfying
$\tau_{ki}<a_{ki}$ such that (\ref{3x1}) holds for all $u\in C_0^\infty
({\H}^n)$. For $\epsilon>0$, we let $\psi_{k, \epsilon} (r)$ be as in Theorem 1.1, that is
 $$ \psi_{k,\epsilon}(r)=\left[ \frac{1+\epsilon^2}{2\epsilon} \big( \cosh \, r- \frac{1-\epsilon^2}
{1+\epsilon^2}\big)\right]^{k-\frac{n}{2}},$$
 where $r$ is the distance to the origin on ${\H}^n$.
  By (\ref{2?1}) and (\ref{2--9}), it follows that
 \begin{eqnarray*}  &  \int_{{\H}^n} |\psi_{k,\epsilon} (r)|^q dV_h =
 \int_{B_n} |G_{k,\epsilon} (x)|^q dx
 =\int_{B_n}
\epsilon^{-n}
\bigg(\frac{1+ \left(\frac{|x|}{\epsilon}\right)^2}{2} \bigg)^{-n} dx\\
  = & \omega_{n-1}
 \int_0^1 \epsilon^{-n}
 \bigg(\frac{1+ \left(\frac{s}{\epsilon}\right)^2}{2}\bigg)^{-n}
  s^{n-1} ds= 2^n \omega_{n-1} \int_0^{\frac{1}{\epsilon}}
    (1+z ^2)^{-n} z^{n-1} dz,\end{eqnarray*}
  where  $$G_{k,\epsilon}
(x)=\bigg[\frac{2\epsilon}{\epsilon^2+|x|^2}\bigg]^{\frac{n}{2}-k}.$$
    Hence  \begin{eqnarray} \label{3-11} \lim_{\epsilon\to 0^+} \int_{{\H}^n}
     |\psi_{k,\epsilon} (r)|^q dV_h &=&
   2^n \omega_{n-1} \int_0^{+\infty} (1+z ^2)^{-n} z^{n-1} dz\\
   &=& 2^{n-1}\omega_{n-1} \frac{\Gamma (n/2)\Gamma(n/2)}{\Gamma(n)}=
 \omega_n,  \nonumber \end{eqnarray}
which also implies that
 $ \int_{{\H}^n} |\psi_{k,\epsilon}(r)|^q dV_h$ is
  decreasing  with respect to $\epsilon$, $\,(\epsilon>0)$.
 According to (\ref{2?4}),
 we have
      $$P_k  \left(\psi_{k,\epsilon}(r)\right) =
\frac{1}{\Lambda_k \omega_n^{2k/n}} \left(\psi_{k,\epsilon}(r)\right)^{q-1},
\quad \,  0\le r<+\infty,$$
 i.e.,
 $$\triangle^{k}_h \psi_{k,\epsilon} (r)+
 \sum_{m=0}^{k-1} a_{km}\triangle_h^{m} \psi_{k,\epsilon} (r)=
 \frac{1}{\Lambda_k \omega_n^{2k/n}} {\big(\psi_{k,\epsilon} (r)\big)}^{\frac{n+2k}{n-2k}},\quad \; 0\le r<+\infty, $$
 where $P_k$ is  the $2k$-order GJMS operator on ${\H}^n$ (see (\ref{2?5}).
 It follows that for any $\epsilon>0$,  \begin{eqnarray*} &&
\frac{\int_{{\H}^n} \big(P_k \psi_{k,\epsilon}(r) \big)
\big(\psi_{k,\epsilon} (r)\big)
 dV_h }{\| \psi_{k,\epsilon} (r)\|_{L^q ({\H}^n)}^2}
   =  \frac{\big(\Lambda_k \omega_n^{2k/n}\big)^{-1} \int_{{\H}^n} \big(
  \psi_{k,\epsilon}(r))^q
 dV_h }{ \big(\int_{{\H}^n} ( \psi_{k,\epsilon} (r))^q
 \big)^{2/q}}\\
&&\quad \;\;  =  \frac{1}{\Lambda_k\omega_n^{2k/n}} \bigg( \int_{{\H}^n}
\big(\psi_{k,\epsilon} (r)\big)^q dV_h \bigg)^{1-2/q}<
  \frac{\omega_n^{1-2/q}}{\Lambda_k \omega_n^{2k/n}} = \frac{1}{\Lambda_k}.\end{eqnarray*}
 Since  $$ T_i \psi_{k,\epsilon}(r) = P_k \psi_{k,\epsilon}(r) +\sum_{m=0}^i
 (\tau_{km}-a_{km}) \triangle_h^m (\psi_{k,\epsilon}(r)),$$
 we obtain that
 for any $\epsilon>0$,
   \begin{eqnarray} \label {3.12}  && \frac{\int_{{\H}^n}
   \big(T_i \psi_{k,\epsilon}(r)\big)
  \big(\psi_{k,\epsilon}(r)\big) dV_h }
 {\| \psi_{k,\epsilon}(r)\|_{L^q({\H}^n)}^2}
   =  \frac{\int_{{\H}^n} \big(P_k
 \psi_{k,\epsilon} (r)\big) \big(\psi_{k,\epsilon}(r)\big) dV_h
}{\| \psi_{k,\epsilon}(r)\|_{L^q({\H}^n)}^2} \\
 && \; \;\quad \quad \quad \; + \frac{\int_{{\H}^n} \bigg(\sum_{m=0}^i (\tau_{km} -
a_{km})\triangle_h^m (\psi_{k,\epsilon}(r))\bigg)
(\psi_{k,\epsilon}(r))dV_h}
{\| \psi_{k,\epsilon}(r) \|_{L^q ({\H}^n)}^2}\nonumber \\
 && \quad \quad \quad  <   \frac{1}{\Lambda_k} +
 \frac{\int_{{\H}^n} \bigg(\sum_{m=0}^i (\tau_{km}
- a_{km})\triangle_h^m (\psi_{k,\epsilon}(r))\bigg)
( \psi_{k,\epsilon}(r))dV_h} {\|\psi_{k,\epsilon} (r)
\|_{L^q ({\H}^n)}^2}. \nonumber \end{eqnarray}
    Clearly, if $i=0$, then (\ref{3.12}) implies
 $$\frac{\int_{{\H}^n} \big(T_i (\psi_{k,\epsilon}(r))\big)
\big(\psi_{k,\epsilon} (r)\big) \, dV_h}
 {\|\psi_{k,\epsilon} (r)\|_{L^{q}({\H}^n)}^2} <\frac{1}{\Lambda_k}.$$
 This contradicts (3.19) because we can choose a function
$v_\epsilon \in C_0^\infty({\H}^n)$
 (see later) such that
$$\frac{\int_{{\H}^n} \big(T_i v_\epsilon\big)
 v_\epsilon \, dV_h}
 {\|v_\epsilon\|_{L^{q}({\H}^n)}^2} <\frac{1}{\Lambda_k}.$$
Thus the conclusion of the theorem holds for $i=0$.

Now, we consider the case $i\ge  1$. It is obvious that
$$ \lim_{r\to +\infty} \frac{\partial^m \psi_{k,\epsilon}(r)}{\partial r^m}=0,\quad
\; \lim_{\underset {r\ne 0} {\epsilon\to 0^+}} \frac{\partial^m
\psi_{k,\epsilon} (r)}{\partial r^m} =0,
 \quad \; m=0,1,\cdots, i-1.$$
It follows from (\ref{3-11}) that  $$ \frac{1}{\|
\psi_{k,\epsilon} (r) \|_{L^q ({\H}^n)}^2}=\frac{1}{\omega_n^{2/q}}+o(\epsilon) \quad \; \mbox{as}\;\;
\epsilon\to 0^{+}.$$ Therefore, for $\epsilon$ sufficiently close to $0$,
 we have
   \begin{eqnarray*}  && \frac{\int_{{\H}^n} \bigg(\sum_{m=0}^i (\tau_{km}
- a_{km})\triangle_h^m (\psi_{k,\epsilon} (r))\bigg)
(\psi_{k,\epsilon}(r))dV_h} {\| \psi_{k,\epsilon} (r)
\|_{L^q
({\H}^n)}^2} \\
&& \; \quad \; \quad  \;\; = \omega_n^{-2/q} \sum_{m=0}^i (\tau_{km} -
a_{km})\int_{{\H}^n} |\triangle_h^{m/2} (\psi_{k,\epsilon}(r))|^2 dV_h +
 f_1 (\epsilon),\end{eqnarray*}
 where $$\lim_{\epsilon\to 0^{+}}
\frac{f_1(\epsilon)}{ \sum_{m=0}^i (\tau_{km} - a_{km})\int_{{\H}^n}
|\triangle_h^{m/2} (\psi_{k,\epsilon}(r))|^2 dV_h}=0.$$ For
each $m\in \{0, 1, \cdots, i-1\}$, it follows from Nirenberg's lemma
(see \cite{HL},  \cite{N} or [1, Lemma 14.1]) that
 for any $\varrho>0$, there exists a constant $c_m$ depending only on
 $\varrho$, $m$ and $i$ such that
\begin{eqnarray}\label {3.17} \int_{{\H}^n} |\triangle_h^{m/2} u|^2 dV_h \le \varrho
\int_{{\H}^n} |\triangle_h^{i/2} u|^2 dV_h
 +c_m \int_{{\H}^n} |u|^2 dV_h\end{eqnarray}
 for all $u\in W^{k,2}({\H}^n)$. We can choose $\varrho>0$
small enough such that
  \begin{eqnarray}\label{3'13} \tau_{ki}-a_{ki} + \varrho \sum_{m=1}^{i-1} |\tau_{km}-a_{km}|<0.\end{eqnarray}
 Thus
   \begin{eqnarray} \label{3c1}  && \quad  \frac{\int_{{\H}^n} \big(T_i (\psi_{k,\epsilon}
 (r)\big)\big(\psi_{k,\epsilon} (r)\big)dV_h }
  {\|\psi_{k,\epsilon} (r)\|_{L^q({\H}^n)}^2} \\
  && \quad \quad \; \;  <   \frac{1}{\Lambda_k} + \omega_n^{-2/q}
 \left[\bigg(\tau_{ki}-a_{ki}+ \varrho \sum_{m=1}^{i-1}
|\tau_{km}-a_{km}|\bigg)
 \int_{{\H}^n} |\triangle^{i/2}_h (\psi_{k,\epsilon}(r))|^2dV_h
 +f_2 (\epsilon) \right. \nonumber \\
  && \left. \quad\;\;\quad \;\;\;
+ c'_i \int_{{\H}^n}|\psi_{k,\epsilon} (r)|^2
dV_h\right],\nonumber \end{eqnarray}
 where $$\lim_{\epsilon\to 0^+} \frac{f_2(\epsilon)}{
\bigg(\tau_{ki}-a_{ki}+ \varrho \sum_{m=1}^{i-1}
|\tau_{km}-a_{km}|\bigg)
 \int_{{\H}^n} |\triangle^{i/2}_h (\psi_{k,\epsilon}(r))|^2dV_h}=0,$$
and $c'_i$ is a constant depending only on $\varrho$ and $i$.

From (\ref{2--9}), we have
   \begin{eqnarray*}  \int_{{\H}^n} |\psi_{k,\epsilon} (r)|^2 dV_h  &=&
\int_{B_n} \left(\frac{2\epsilon}{1+\epsilon^2}\right)^{n-2k}
 \bigg[ \frac{2\epsilon^2}{1+\epsilon^2}
+\frac{2|x|^2}{1-|x|^2}\bigg]^{2k-n}
\bigg(\frac{2}{1-|x|^2}\bigg)^n  dx\\
 &=& 2^n\omega_{n-1}\left(\frac{2\epsilon}{1+\epsilon^2}\right)^{n-2k}
 \int_0^{1}
 \left[ \frac{2\epsilon^2}{1+\epsilon^2}  (1-s^2) +2s^2\right]^{2k-n}
 \left(\frac{1}{1-s^2}\right)^{2k} s^{n-1} ds.\end{eqnarray*}
   Making the change of variables $t=\frac{s}{\epsilon}$,
 we obtain that
\begin{eqnarray*}  \int_{{\H}^n} |\psi_{k,\epsilon} (r)|^2 dV_h &=& 2^n\omega_{n-1} \left(\frac{2\epsilon}{1+\epsilon^2}\right)^{n-2k}
\int_0^{1/\epsilon}
  \big[\frac{2\epsilon^2}{1+\epsilon^2}  (1+t^2)\big]^{2k-n} \bigg(\frac{1}{1- \epsilon^2 t^2 }\bigg)^{2k}
   \epsilon^n t^{n-1} dt\\
   &=& 2^n\omega_{n-1}  \epsilon^{2k} \left(\frac{1+\epsilon^2}{2\epsilon} \right)^{2}
    \int_0^{1/\epsilon}
   \big(\frac{\big(\frac{2\epsilon}{1+\epsilon^2}\big)^2}{1-\epsilon^2 t^2}\big)\big(\frac{1}{1-
\epsilon^2 t^2 }\big)^{2k-1}
 \frac{t^{n-1} dt}{\big(1+t^2\big)^{n-2k}}\end{eqnarray*}
  It follows from (\ref{2--9}) that $$\frac{\partial \psi_{k,\epsilon}(r)}{\partial x_i} =
  \left(\frac{2\epsilon}{1+\epsilon^2}\right)^{\frac{n}{2}-k}
 \big(k-\frac{n}{2}\big)
\left[\frac{2\epsilon^2}{1+\epsilon^2}  +\frac{2|x|^2}{1-|x|^2}\right]^{k-\frac{n}{2}-1} \frac{4x_i}{(1-|x|^2)^2},$$
    which implies
      \begin{eqnarray*} &&\int_{{\H}^n} |\nabla_{h} \psi_{k,\epsilon}(r)|^2 dV_{h} \\
  &=&  \left(\frac{2\epsilon}{1+\epsilon^2}\right)^{n-2k}  \big(k-\frac{n}{2}\big)^2 \int_{B_n}
\left[ \frac{2\epsilon^2}{1+\epsilon^2}    +\frac{2|x|^2}{1-|x|^2}\right]^{2k-n-2} |x|^2 \big(\frac{2}{1-|x|^2}\big)^{n+2} dx\\
 &=& 2^{n+2}  \left(\frac{2\epsilon}{1+\epsilon^2}\right)^{n-2k}  \big(k-\frac{n}{2}\big)^2 \omega_{n-1} \\
 &&  \;\; \quad \;\; \quad \times  \int_0^{1} \left[ \frac{2\epsilon^2}{1+\epsilon^2}   (1-s^2)+2s^2\right]^{2k-n-2}
\big(\frac{1}{1-s^2}\big)^{2k} s^{n+1}ds. \end{eqnarray*}
   Again, by the substitution
 $t=\frac{s}{\epsilon}$, we have
   \begin{eqnarray*} \int_{{\H}^n} |\nabla_{h} \psi_{k,\epsilon} (r)|^2 dV_{h}&=
 &  2^{n+2}  \big(k-\frac{n}{2}\big)^2 \omega_{n-1}  \epsilon^{2k}
 \left(\frac{1+\epsilon^2}{2\epsilon} \right)^{2}\qquad \quad\qquad
\\
  &&   \times \int_0^{\frac{1}{\epsilon}}
   \bigg(\frac{1}{1-\epsilon^2 t^2} \bigg)^{2k}
 \frac{t^{n+1}dt}{\big(1+t^2\big)^{n-2k+2}}.\qquad \quad\qquad\end{eqnarray*}
   It follows from the dominated convergence theorem, if $n>4k-2$,
 $$ \lim_{\epsilon\to 0^+}\int_0^{\frac{1}{\epsilon}}
 \bigg(\frac{\left(\frac{2\epsilon^2}{1+\epsilon^2}\right)^2}{1-\epsilon^2t^2}\bigg)
  \bigg(\frac{1}{1- \epsilon^2 t^2 }\bigg)^{2k-1}
 \frac{t^{n-1}dt}{(1+t^2)^{n-2k}} =0,$$
 and
\begin{eqnarray*} \lim_{\epsilon\to 0^+} \int_0^{1/\epsilon}  \bigg(\frac{1}{1-\epsilon^2  t^2}\bigg)^{2k}
  \frac{t^{n+1} dt}{\big(1+t^2\big)^{n-2k+2}}=\int_0^{+\infty}
 \frac{t^{n+1} dt} {\big(1+t^2)^{n-2k+2}}= \frac{\Gamma(\frac{n}{2}+1) \Gamma(\frac{n}{2} +1-2k)}{2\,\Gamma(n+2-2k)}.\end{eqnarray*}
The latter integral on the right-hand side is a finite positive constant for $n>4k-2$.
 On the other hand,  since
  $$\lim_{\epsilon\to 0^+} \frac{2^n\omega_{n-1}
 \epsilon^{2k} \left(\frac{1+\epsilon^2}{2\epsilon} \right)^{2}}
{2^{n+2}  \big(k-\frac{n}{2}\big)^2 \omega_{n-1}
 \epsilon^{2k} \left(\frac{1+\epsilon^2}{2\epsilon} \right)^{2}}=\frac{1}{4(k-\frac{n}{2})^2},$$
  it follows that
     \begin{eqnarray}\label{3-20} \lim_{\epsilon\to 0^+}\frac{\int_{{\H}^n} |\psi_{k,\epsilon}(r)|^2 dV_h}
    {\int_{{\H}^n}
 |\nabla_h (\psi_{k,\epsilon}(r))|^2 dV_h}=0.\end{eqnarray}
 By taking $m=1$ and replacing $u$ by $\psi_{k,\epsilon}(r)$ in
(\ref{3.17}), we have
 $$ 1\le   \varrho \frac{\int_{{\H}^n} |\triangle_h^{i/2}
 (\psi_{k,\epsilon} (r))|^2 dV_h}{\int_{{\H}^n}
|\nabla_h ( \psi_{k,\epsilon} (r))|^2 dV_h}
 +c_1 \frac{\int_{{\H}^n} |\psi_{k,\epsilon}(r)|^2 dV_h}{\int_{{\H}^n}
|\nabla_h (\psi_{k,\epsilon}(r))|^2 dV_h}.$$
 Letting  $\epsilon \to 0^+$, we find by this and (\ref{3-20}) that
 \begin{eqnarray}\label{3-21}  \underset {\epsilon\to 0^+}  {\underline{\lim}} \frac{\int_{{\H}^n}
 |\triangle_h^{i/2} \big(\psi_{k,\epsilon}(r)\big)|^2 dV_h}
{\int_{{\H}^n} |\nabla_h \big(\psi_{k,\epsilon}(r)\big)|^2
dV_h} \ge \frac{1}{\varrho}.\end{eqnarray}
   Combining (\ref{3.17}), (\ref{3'13}), (\ref{3-20}) and (\ref{3-21}),
   we obtain that if $\tau_{ki}<a_{ki}$ and $n>4k-2$,
then for $\epsilon$ sufficiently close to $0$,
$$ \left(\tau_{ki}-a_{ki} + \varrho  \sum_{m=0}^{i-1}
|\tau_{km}-a_{km}| \right) \int_{{\H}^n}
|\triangle_h^{i/2}\big(\psi_{k,\epsilon}(r)\big)|^2 dV_h +
 c'_i \int_{{\H}^n}|\psi_{k,\epsilon}(r)|^2 dV_h <0.$$
 From (\ref{3c1}), it follows that
 for $\epsilon>0$ sufficiently close to $0$,
 $$ \int_{{\H}^n} \big(T_i (\psi_{k,\epsilon}(r))\big)
 \big(\psi_{k,\epsilon}(r)\big)dV_h  <  \frac{1}{\Lambda_k}
 \|\psi_{k,\epsilon}(r)\|_{L^q({\H}^n)}^2.$$

 Finally, since $G_{k,\epsilon} \in W^{k,2}({\Bbb R}^n)$ (see, Theorem 1.1 of \cite{CT}, or
 \cite{B}), by (\ref{2--9}) we get that $\psi_{k,\epsilon} (y)\in
 W^{k,2}({\Bbb H}^n)$, $\, \epsilon>0$. Thus,
  for $\epsilon$ sufficiently close to $0$,
   we can choose a sequence of smooth functions $w_j$ with compact support in
${\H}^n$ such that
$$\lim_{j\to +\infty} \|\psi_{k,\epsilon}-w_j\|_{W^{k,2}({\H}^n)}=0.$$
 It follows that there is an integer $j_0$ such that
 $$\int_{{\H}^n} \big(T_i (w_{j_0})\big)
 \big(w_{j_0}\big)dV_h < \frac{1}{\Lambda_k} \|w_{j_0}\|_{L^q({\H}^n)}^2.$$
  Let us denote $w_{j_0}$ by $u$.
 This is in contradiction with (\ref{3x1}), which proves the theorem when $n>4k-2$.$\;\;\square$

\vskip 1.89  true cm

{\bf Acknowledgments}

 \vskip 0.58 true cm

I wish to express my sincere gratitude to Professor
 L. Nirenberg, Professor Fang-Hua Lin and Professor YanYan Li
 for their interest and support in this project.
  I would also like to thank Professor T. Aubin for some useful
  comments and for pointing out reference \cite{CT} to me.
  Finally, I owe special thanks to the referee,
 whose many helpful comments and constructive suggestions have greatly
 improved the quality of this paper. This
research was supported by SRF for ROCS, SEM (No. 2004307D01) and NNSF of China
(No: 11171023/A010801).

\vskip 2.0 true cm

\vskip 0.32 true cm

\end{document}